# H-Matrix Accelerated Direct Matrix Solver for Maxwell's Equations using the Chebyshev-based Nyström Boundary Integral Equation Method

JIN HU[1](Student Member, IEEE), EMRAH SEVER[2], OMID BABAZADEH[3](Student Member, IEEE), IAN JEFFREY[3](Member, IEEE), VLADIMIR OKHMATOVSKI[3](Senior Member, IEEE), CONSTANTINE SIDERIS[1](Senior Member, IEEE)

[1] Department of Electrical and Computer Engineering, University of Southern California, USA
[2] ASELSAN Inc., Ankara, Turkey
[3] Department of Electrical and Computer Engineering, University of Manitoba, Canada

CORRESPONDING AUTHOR: C. Sideris (e-mail: csideris@usc.edu).

**ABSTRACT** An H-matrix accelerated direct solver employing the high-order Chebyshev-based Boundary Integral Equation (CBIE) method has been formulated, tested, and profiled for performance on high contrast dielectric materials and electrically large perfect electric conductor objects. The matrix fill performance of the CBIE proves to be fast for small to moderately sized problems compared to its counterparts, e.g. the locally corrected Nyström (LCN) method, due to the way it handles the singularities by means of a global change of variable method. However, in the case of electrically large scattering problems, the matrix fill and factorization still dominate the solution time when using a direct solution approach. To address this issue, an H-Matrix framework is employed, effectively resolving the challenge and establishing the CBIE as a competitive high-order method for solving scattering problems with poorly conditioned matrix equations. The efficacy of this approach is demonstrated through extensive numerical results, showcasing its robustness to problems that are electrically large, near physical resonances, or that have large dielectric permittivities. The capability of the proposed solver for handling arbitrary geometries is also demonstrated by considering various scattering examples from complex CAD models.

## I. INTRODUCTION

Boundary integral equations (BIE) are the most common methods for solving large-scale electromagnetic scattering problems in practice since they reduce the dimension of the problem by one and only require the discretization of the contour or surface for 2D and 3D problems, respectively, making them more favorable compared to their volumetric counterparts. In addition, BIEs satisfy the radiation condition inherently as the appropriate Green's function is the kernel of the integral operator. Among the methods for solving BIEs, formulations using the low-order method of moments (MoM), i.e., flat triangular patches and low-order RWG basis functions, are widely used due to the relative ease of implementation and efficiency for producing a discretization of the BIE. Despite low-order approaches being typical for large electromagnetic scattering problems, high-order methods continue to attract interest because i) they reduce the number of unknowns ii) they offer better representation of the geometry through high-order curved patches, and iii) they enable a path towards controlling the error in the solution. As industries continue to push the complexity and electrical sizes in their designs, high-order methods become critical tools for accurate simulation.

The Nyström method (NM) [1] is one of several high-order frameworks for solving BIEs. However, for the singular integral kernels arising in electromagnetics, a standard NM approach fails to accurately capture the physics. To overcome this problem, the locally corrected Nyström (LCN) method was proposed in [2], where a local correction to the integration weights is introduced to accurately capture the behavior of the kernel near singularities. The details of a high-order formulation of this method were presented in [3], where it was applied to the electric field integral equation (EFIE) and magnetic field integral equation (MFIE). Therein, and in [4], the equivalence of LCN to a high-order implementation of MoM was shown. From an implementation perspective, high-order MoM requires the evaluation of integrals over both the expansion and testing functions for matrix fill, incurring additional computational cost compared to LCN. However, since the standard basis formulation used in LCN does not enforce a continuity condition across element edges, a combination of both methods might be beneficial [4][5][6][7].

Recently, a new NM formulation referred to as the Chebyshev-based Boundary Integral Equation (CBIE) was proposed in [8] for acoustic problems and applied to electromagnetic scattering problems in [9]. In this NM approach, the surface current is expanded in terms of Chebyshev polynomials and the singularities of the kernels of the integral equation are handled by a singularity-cancelling change of variables. Compared to the LCN method, the CBIE is significantly faster when filling the matrix since the LCN method requires the solution of a set of local correction systems accounting for singularities in nearby interactions between observation points and source patches. The size of these local systems is proportional to the square of the order of quadrature method, thereby increasing the numerical burden as the order of the method increases. In addition, adaptive integration is required for the singular and nearly singular interactions to control the level of accuracy, a process that dominates the computation time as the order of the quadrature and/or the number of patches increases. These drawbacks cause LCN to be slow when the electrical size of the problems increases compared to the wavelength. In contrast, the CBIE's singularity cancelling change of variable method for handling singularities can be more efficient, resulting in fast matrix fills for high-order solutions [9].

Both the LCN and CBIE methods require the solution of an algebraic system of equations. The solution of the system can be obtained from a direct solver (i.e., LU decomposition) or an iterative method like GMRES. While iterative solvers may be preferred (or deemed necessary) for large systems of equations, their convergence generally depends on the condition number of the system. In cases of a frequency near a physical resonance or when geometric singularities are featured in the scatterer, iterative solvers often fail to solve the system accurately. Thus, if direct solvers become necessary for a particular problem, the difficulty shifts to ensuring that the direct solver framework is efficient, avoiding the need for $O(N^3)$ and $O(N^2)$ compute time and storage for a system involving N unknowns. While the CBIE accelerates matrix fill times, filling the matrix entirely still requires $O(N^2)$ operations and memory and direct factorization of the matrix requires $O(N^3)$ operations, which make it prohibitive for large scale problems.

The goal of our current work is to combine CBIE with a hierarchical matrix (H-Matrix) framework to produce an efficient direct solver for large-scale scattering problems. This paper extends the initial work presented in [10], in which we demonstrated preliminary results using the MFIE for PEC objects. H-matrix acceleration of IE solutions also has been previously demonstrated for the single source surface-volume-surface (SS-SVS) integral equation for scattering by dielectric objects [11], composite dielectric objects in free space [12], in layered media [13] and composite metal-dielectric objects [14] discretized by RWG based method of moments. The application of the H-matrix framework to the LCN discretization of the MFIE was presented in [15]. Extending the work presented in [10], the H-matrix method and implementation approach is described in detail and is also applied to the N-Müller integral equation of scattering by a homogeneous dielectric object filled with high contrast material discretized by the CBIE method. Specifically, CBIE provides fast and efficient computation of matrix entries, easily scalable to high-order solutions, while H-matrices provide an efficient path to accelerate the direct solution of the resulting algebraic system.

The organization of the paper is as follows. In Section II the formulation of the CBIE for the MFIE is given, while in Section III the theoretical details of the application of H-Matrices to the CBIE are provided. Section IV presents numerical results demonstrating the benefits of the CBIE/H-Matrix framework, and the paper is concluded in Section V. An $e^{j\omega t}$ time dependency is assumed and suppressed throughout the paper.

## II. DISCRETIZATION OF THE INTEGRAL EQUATIONS WITH CBIE METHOD

### A. INTEGRAL EQUATIONS

Starting from Green's identities and using the equivalence principle one can derive the following integral representations [16] in terms of the equivalent electric surface current density **J** and magnetic surface current density **M** for the outer and inner domains of a homogeneous penetrable body $\Omega_2$ surrounded by the surface $S$ immersed in homogeneous domain $\Omega_1$ as given in Fig. 1.

$$
\begin{aligned}
\mathbf{E}_1(\mathbf{r}) &= \mathbf{E}^{\text{inc}}(\mathbf{r}) - L_E^+\{\mathbf{J},\mathbf{M}\}(\mathbf{r}), &\mathbf{r} \in \Omega_1\backslash S \\
\mathbf{H}_1(\mathbf{r}) &= \mathbf{H}^{\text{inc}}(\mathbf{r}) - L_H^+\{\mathbf{J},\mathbf{M}\}(\mathbf{r}), &\mathbf{r} \in \Omega_1\backslash S \\
\mathbf{E}_2(\mathbf{r}) &= L_E^-\{\mathbf{J},\mathbf{M}\}(\mathbf{r}), &\mathbf{r} \in \Omega_2\backslash S \\
\mathbf{H}_2(\mathbf{r}) &= L_H^-\{\mathbf{J},\mathbf{M}\}(\mathbf{r}), &\mathbf{r} \in \Omega_2\backslash S
\end{aligned}
\quad (1)
$$

These integral representations are also known as Stratton-Chu formulas [17]. The $L_E$ and $L_H$ operators in (1) are defined as follows:

$$
\begin{aligned}
L_E^\pm\{\mathbf{J},\mathbf{M}\} &= j\omega\mu_\pm \mathcal{S}_\pm\{\mathbf{J}\} - \frac{1}{j\omega\epsilon_\pm}\nabla\mathcal{S}_\pm\{\nabla'_s \cdot \mathbf{J}\} \\
&\quad + \nabla \times \mathcal{S}_\pm\{\mathbf{M}\} \\
L_H^\pm\{\mathbf{J},\mathbf{M}\} &= j\omega\epsilon_\pm \mathcal{S}_\pm\{\mathbf{M}\} - \frac{1}{j\omega\mu_\pm}\nabla\mathcal{S}_\pm\{\nabla'_s \cdot \mathbf{M}\} \\
&\quad - \nabla \times \mathcal{S}_\pm\{\mathbf{J}\}
\end{aligned}
\quad (2)
$$

where $\mathcal{S}_\pm$ is the well-known single layer potential and is defined for a scalar or vectoral function $\sigma$ as

$$
\mathcal{S}_\pm[\sigma](\mathbf{r}) = \int_S G_\pm(\mathbf{r},\mathbf{r}')\sigma(\mathbf{r}')dS' \quad (3)
$$

where

$$G_{\pm}(\mathbf{r}, \mathbf{r}') = \frac{e^{-jk_{\pm}|\mathbf{r}-\mathbf{r}'|}}{4\pi|\mathbf{r}-\mathbf{r}'|} \quad (4)$$

is the 3-dimensional free space Green's function.

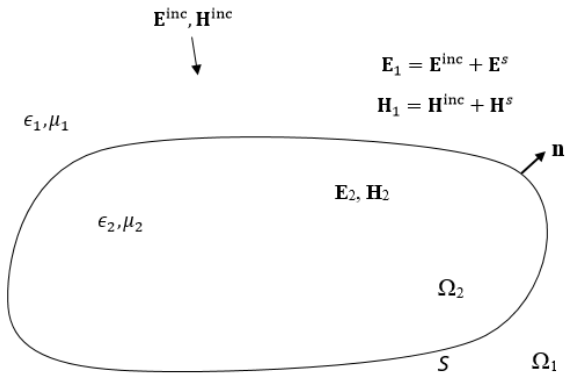

**FIGURE 1.** The geometry of the boundary value problem.

In equations (1)-(4) the $\pm$ indicates approaching the boundary S from $\Omega_1$ and $\Omega_2$, respectively. The curly braces used in operators $L^{\pm}\{\mathbf{J}, \mathbf{M}\}$ in (2) indicate that the observation point is not located on the boundary since the third term, $\nabla \times \mathcal{S}_{\pm}$, is not continuous across the boundary and must be defined in a Cauchy principal value sense. Note that square brackets ($[\cdot]$) will be used with continuous operators where the observation point is located on the surface. When the observation point approaches the boundary $S$ (in the limit sense), the term, $\nabla \times \mathcal{S}_{\pm}$, is expressed as a jump condition [18]:

$$\lim_{\substack{h \to 0 \\ h > 0}} \nabla \times \mathcal{S}\{\mathbf{a}\}(\mathbf{r_0} \pm \mathbf{n}h)$$
$$= \mp \frac{1}{2} \mathbf{n} \times \mathbf{a} + \nabla \times \mathcal{S}[\mathbf{a}](\mathbf{r_0}), \quad (5)$$
$$\mathbf{r_0} \in S$$

where $\mathbf{n}$ is the unit outward normal. Due to the relation given by (5), and considering that $\mathbf{n} \times \mathbf{n} \times \mathbf{a} = -\mathbf{a}$, the operators in (2) also have the jump terms when $\mathbf{r} = \mathbf{r_0} \in S$ which are defined as

$$\lim_{\substack{h \to 0 \\ h > 0}} \mathbf{n} \times L_{\mathbf{E}}\{\mathbf{J}, \mathbf{M}\}(\mathbf{r_0} \pm \mathbf{n}h)$$
$$= \pm \frac{\mathbf{M}}{2} + \mathbf{n} \times L_{\mathbf{E}}^{\pm}[\mathbf{J}, \mathbf{M}](\mathbf{r_0}),$$
$$\mathbf{r_0} \in S$$
$$\lim_{\substack{h \to 0 \\ h > 0}} \mathbf{n} \times L_{\mathbf{H}}\{\mathbf{J}, \mathbf{M}\}(\mathbf{r_0} \pm \mathbf{n}h)$$
$$= \mp \frac{\mathbf{J}}{2} + \mathbf{n} \times L_{\mathbf{H}}^{\pm}[\mathbf{J}, \mathbf{M}](\mathbf{r_0}),$$
$$\mathbf{r_0} \in S$$
(6)

On a penetrable boundary, the relation between surface currents densities and fields is given by $\mathbf{J} = \mathbf{n} \times \mathbf{H_1} = \mathbf{n} \times \mathbf{H_2}$ and $\mathbf{M} = \mathbf{E_1} \times \mathbf{n} = -\mathbf{n} \times \mathbf{E_2}$ where $\mathbf{n}$ is the unit surface outward normal directed into $\Omega_1$ as shown in Figure 1.

If the equations in (1) are crossed with $\mathbf{n}$ to get the tangential components of the fields and the above-given relations are imposed as boundary conditions, then the EFIE and MFIE of a penetrable body are obtained as

$$L_{\mathbf{E}_{tan}}^{\pm}[\mathbf{J}, \mathbf{M}] = \mp \frac{\mathbf{M}}{2} + \mathbf{n} \times L_{\mathbf{E}}^{\pm}[\mathbf{J}, \mathbf{M}]$$
$$- \begin{cases} \mathbf{n} \times \mathbf{E}^{inc} \\ 0 \end{cases} = 0,$$
$$L_{\mathbf{H}_{tan}}^{\pm}[\mathbf{J}, \mathbf{M}] = \pm \frac{\mathbf{J}}{2} + \mathbf{n} \times L_{\mathbf{H}}^{\pm}[\mathbf{J}, \mathbf{M}]$$
$$- \begin{cases} \mathbf{n} \times \mathbf{H}^{inc} \\ 0 \end{cases} = 0.$$
(7)

where the subscript "*tan*" indicates that the tangential components of the fields are taken. In (7) there are two unknowns in a system of four vector equations. To solve the system, any two of these four equations can be chosen or a linear combination of them can be set as

$$L_{\mathbf{E}_{tan}}^{+}[\mathbf{J}, \mathbf{M}](\mathbf{r_0}) + \alpha L_{\mathbf{E}_{tan}}^{-}[\mathbf{J}, \mathbf{M}](\mathbf{r_0}) = 0,$$
$$L_{\mathbf{H}_{tan}}^{+}[\mathbf{J}, \mathbf{M}](\mathbf{r_0}) + \beta L_{\mathbf{H}_{tan}}^{-}[\mathbf{J}, \mathbf{M}](\mathbf{r_0}) = 0.$$
(8)

Each integral equation system will have different spectral properties depending on which kind of kernels they have. Choosing $\alpha = \beta = 1$ results in the well-known PMCHWT formulation for a penetrable object. Alternatively, the choice

$$\alpha = -\frac{\epsilon_2}{\epsilon_1}, \quad \beta = -\frac{\mu_2}{\mu_1} \quad (9)$$

results in another well-known N-Müller formulation for a penetrable object [16]. If the boundary is a perfect electric conductor (PEC), then $\mathbf{M} = \mathbf{0}$ in (5) and using (2) we get the following classical electric and magnetic field integral equations for a PEC object:

$$L_{\mathbf{E}_{tan}}^{+}[\mathbf{J}, \mathbf{0}](\mathbf{r_0}) = \mathbf{n} \times \left[ j\omega\mu_1 \mathcal{S}_{+}[\mathbf{J}](\mathbf{r_0}) \right.$$
$$\left. - \frac{1}{j\omega\epsilon_1} \nabla \mathcal{S}_{+}[\nabla' \cdot \mathbf{J}](\mathbf{r_0}) \right]$$
$$- \mathbf{n} \times \mathbf{E}^{inc}(\mathbf{r_0}) = 0,$$
$$L_{\mathbf{H}_{tan}}^{+}[\mathbf{J}, \mathbf{0}](\mathbf{r_0}) = \frac{1}{2}\mathbf{J}(\mathbf{r_0}) + \mathbf{n} \times [\nabla \times \mathcal{S}_{+}[\mathbf{J}](\mathbf{r_0})]$$
$$- \mathbf{n} \times \mathbf{H}^{inc}(\mathbf{r_0}) = 0,$$
$$\mathbf{r_0} \in S.$$
(10)

The EFIE and MFIE for both penetrable and PEC objects given by (7) and (10) may include resonance solutions associated with the inner boundary value problem arising from the equivalence principle. However, the PMCHWT and Müller integral equations, which are the combination of the integral equations of the inner and outer domains, are resonance free. The PMCHWT is a Fredholm IE of the first kind and contains a strongly singular kernel. On the other hand, the Müller integral equation is a Fredholm IE of the second kind and less singular since the strong singularity is cancelled due to the proper combination of the integrals by coefficients in (7). In this work we consider applying H-matrix accelerated CBIE formulation to the MFIE and the Müller formulations for scattering from PEC objects and penetrable dielectric objects respectively.

To evaluate the effectiveness of the CBIE accelerated by H-matrices as a direct solver, two significant cases are explored in this study. First, we demonstrate that the proposed method efficiently solves the scattering problem of a sphere filled with high-contrast material exhibiting a resonant behavior. In contrast, an iterative solver struggles and requires a high number of iterations to get comparable accuracy in the solution. Second, a dielectric humanoid bugs bunny model is solved using both the proposed method and an iterative solver, and it is shown that the number of iterations and hence the solution time required by the iterative solver increases linearly as the dielectric permittivity increases. In addition to these dielectric scenarios, the practical utility of the method is showcased by solving the MFIE for several different PEC objects, such as a sphere, a B2-aircraft, and a spiral antenna.

### B. DISCRETIZATION OF THE IEs BY CBIE

As a surface integral equation method, the CBIE first discretizes the surface into a number of patches, on each of which it expands the unknown surface current density in terms of unitary vectors defined in a local $(u, v)$ coordinate system as

$$\boldsymbol{F}^p(u,v) = \frac{F^{p,u}(u,v)\boldsymbol{a}_u^p(u,v) + F^{p,v}(u,v)\boldsymbol{a}_v^p(u,v)}{\sqrt{|G^p(u,v)|}} \quad (11)$$

where $\boldsymbol{F}^p$ represents the surface current $\boldsymbol{J}^p$ or $\boldsymbol{M}^p$ defined in terms of the local coordinate system on the $p$th patch, $\boldsymbol{a}_{u,v}^p$ are the covariant basis vectors, and $\sqrt{|G^p|}$ is the Jacobian of the surface [9]. Each component of the current is expanded in a set of orthogonal polynomials. The CBIE uses Chebyshev polynomials due to their fast-converging behavior and favorable discrete orthogonality properties. The scalar component of the discretized current density on the surface is represented in terms of Chebyshev polynomials ($T_n(\cdot)$) as

$$F^{p,a}(u,v) = \sum_{m=0}^{N_v^p-1} \sum_{n=0}^{N_u^p-1} \gamma_{n,m}^{p,a} T_n(u) T_m(v), \quad (12)$$
$$\text{for } a = u, v$$

where $F^{p,a} (a = u, v)$ are the scalar contravariant components of the surface current, and $\gamma_{n,m}^{p,a}$ are the Chebyshev coefficients which can be computed from the values of the densities on Chebyshev nodes as follows:

$$\gamma_{n,m}^{p,a} = \frac{\alpha_n \alpha_m}{N_u^p N_v^p} \sum_{k=0}^{N_v^p-1} \sum_{l=0}^{N_u^p-1} J^{p,a}(u_l, u_k) T_n(u_l) T_m(v_k). \quad (13)$$

The unitary vectors and the surface normal are defined as

$$\boldsymbol{a}_u^p = \frac{\partial \boldsymbol{r}^p(u,v)}{\partial u}, \quad \boldsymbol{a}_v^p = \frac{\partial \boldsymbol{r}^p(u,v)}{\partial v},$$
$$\hat{\boldsymbol{n}}^p = \frac{\boldsymbol{a}_u^p \times \boldsymbol{a}_v^p}{\|\boldsymbol{a}_u^p \times \boldsymbol{a}_v^p\|}, \quad (14)$$

where $\boldsymbol{r}^p(u,v)$ is the position vector defined on local $(u,v)$ coordinate system.

To evaluate the integrals over the curvilinear patches, Fejer's first quadrature rule, which uses Chebyshev nodes as quadrature points, is used. One important advantage of Chebyshev polynomials is that since they include the endpoints as zeros, they can be used to construct a two-dimensional quadrature rule that includes the contour of the patches which is a requirement when the continuity of the current across the patches is considered. To evaluate the integral over the entire surface we need to discretize the integral operators in (8) and (10) as a sum of integrals over the $M$ patches as

$$\mathcal{K}[\boldsymbol{F}](\boldsymbol{r}) = \sum_{p=1}^{M} \mathcal{K}[\boldsymbol{F}^p](\boldsymbol{r}) \quad (15)$$

where $\boldsymbol{F}^p$ stands for the equivalent electric and magnetic surface current densities. The operator $\mathcal{K}$ may have one of the following forms:

$$\begin{aligned} &\boldsymbol{n} \times \mathcal{S}[\boldsymbol{F}](\boldsymbol{r}_0), \\ &\boldsymbol{n} \times \nabla \times \mathcal{S}[\boldsymbol{F}](\boldsymbol{r}_0), \\ &\boldsymbol{n} \times \nabla \mathcal{S}[\nabla' \cdot \boldsymbol{F}](\boldsymbol{r}_0), \end{aligned} \quad (16)$$

or a difference of them depending on the used integral equation form. The above-given operators have different singularity properties. The first operator is weakly singular, the second is a compact operator, and the third one is a strongly singular kernel which requires more effort to handle the singularity. If one wants to solve the MFIE for a PEC object, then the second equation in (10) is adopted, with the operator in second row of (16) as the kernel. On the other hand, if scattering by a dielectric object is under consideration, then the Müller IE given by (8) and (9) can be chosen. This formulation involves a difference of two operators given by the third row of (16), where one is represented in terms of the parameters of outer domain and the other in terms of the inner domain parameters. The subtraction of these operators results in a less singular kernel due to cancellation of the strongly singular terms [19]. The rest can be handled by a singularity cancellation technique, which, in the CBIE approach, is accomplished by a change of variables [9].

Equations (8) and (10) are vectorial. To produce a system of linear scalar equations, test vectors are used to sample the integral equation. In the CBIE approach, the test vector is chosen as $\sqrt{G^p} \boldsymbol{a}^{p,a} (a = u, v)$ where $\boldsymbol{a}^{p,a}$ are the reciprocal unitary (contravariant) vectors. The contravariant basis vectors $\boldsymbol{a}^{p,u}$ and $\boldsymbol{a}^{p,v}$ are defined via the orthogonality relation [20]

$$\boldsymbol{a}^{p,a} \cdot \boldsymbol{a}_b^p = \begin{cases} 1, & a = b \\ 0, & a \neq b \end{cases} \quad (17)$$

and defined in terms of covariant basis vectors as

$$\boldsymbol{a}^{p,u} = \frac{\boldsymbol{a}_v^p \times \boldsymbol{a}_n^p}{V}, \qquad \boldsymbol{a}^{p,v} = \frac{\boldsymbol{a}_n^p \times \boldsymbol{a}_u^p}{V}, \qquad (18)$$

where $V = \boldsymbol{a}_u^p \cdot (\boldsymbol{a}_v^p \times \boldsymbol{a}_n^p)$ and $\boldsymbol{a}_n^p = (\boldsymbol{a}_u^p \times \boldsymbol{a}_v^p)/V$.

Here, the discretization of the Müller IE (8) will be considered; the MFIE (10) can be discretized by analogy. Once the vector integral equation systems (8) have been discretized as (15) and dotted with $\sqrt{G^p}\boldsymbol{a}^{p,a}$, then on the $p$th patch, we get the algebraic system of the equations in matrix form as

$$\begin{bmatrix} K_{uu}^{J,E} & K_{uv}^{J,E} & K_{uu}^{M,E} & K_{uv}^{M,E} \\ K_{vu}^{J,E} & K_{vv}^{J,E} & K_{vu}^{M,E} & K_{vv}^{M,E} \\ K_{uu}^{J,H} & K_{uv}^{J,H} & K_{uu}^{M,H} & K_{uv}^{M,H} \\ K_{vu}^{J,H} & K_{vv}^{J,H} & K_{vu}^{M,H} & K_{vv}^{M,H} \end{bmatrix} \begin{bmatrix} \mathbf{J}^u \\ \mathbf{J}^v \\ \mathbf{M}^u \\ \mathbf{M}^v \end{bmatrix}$$
$$= \begin{bmatrix} -\epsilon_1 \boldsymbol{a}_v^p \cdot \mathbf{E}^{p,\mathrm{inc}} \\ \epsilon_1 \boldsymbol{a}_u^p \cdot \mathbf{E}^{p,\mathrm{inc}} \\ -\mu_1 \boldsymbol{a}_v^p \cdot \mathbf{H}^{p,\mathrm{inc}} \\ \mu_1 \boldsymbol{a}_u^p \cdot \mathbf{H}^{p,\mathrm{inc}} \end{bmatrix} \qquad (19)$$

where $J^{u,v}$ and $M^{u,v}$ are the contravariant components of the equivalent surface currents and the matrix entries $K_{\cdot\cdot}^{\cdot,E}$ and $K_{\cdot\cdot}^{\cdot,H}$ represents the matrix element sampled from the electric field and magnetic field integral equations, respectively. To construct a high-order-accurate algorithm, each term in (19) should be expanded as

$$K_{ba}^{qp} = K_{ba}^{qp(\mathrm{far})} + K_{ba}^{qp(\mathrm{self})} + K_{ba}^{qp(\mathrm{near})}; \quad b,a = u,v \qquad (20)$$

which involves distinguishing the matrix entries based on the distance between the observation and source patches. If the patches are "far" enough, then the integrals can be calculated accurately by the underlying quadrature method. If the patches coincide (self-patch) or are close to each other (near patches) then a special technique is required to evaluate the integrals since the kernels given by (16) will be singular (or close to singular) and a regular quadrature rule will fail to calculate the integrals with a controlled level of accuracy. The singularity treatment strategy is also different when the observation point sits on the "self" patch or in the "near" vicinity of source patch. These cases can be readily treated by a change of variables method. While details are omitted for brevity, an interested reader is referred to the detailed explanation in [9].

### III. THE H-MATRIX FRAMEWORK

H-matrices leverage a multilevel partitioning of the discretized geometry to categorize the corresponding matrix blocks into low-rank and full-rank sub-matrices where low-rank matrices can be stored in a compressed form. Based on this hierarchical rank decomposition, fast formatted matrix addition and multiplication operations applied to these sub-blocks form the foundation of efficient direct decompositions and solution methods. To determine the low-rank and full-rank sub-blocks, a hierarchical block cluster-tree structure $\mathcal{T}_{I \times I}$ is constructed from the index set $I = \{f_i\}_{i=1,2,\ldots N}$ of the test/basis functions where N is the number of test/basis functions. The procedure starts with the full index set (over all unknowns), which is referred to as the root cluster, and proceeds recursively into sub-blocks [21]. Blocks containing source and observation points that are sufficiently removed in space represent interactions that can be efficiently compressed. These blocks are identified by means of the admissibility condition:

$$\min\{\mathrm{diam}(\mathcal{B}_1), \mathrm{diam}(\mathcal{B}_2)\} \leq \eta \,\mathrm{dist}(\mathcal{B}_1, \mathcal{B}_2) \qquad (21)$$

where $\mathcal{B}_1$ and $\mathcal{B}_2$ are the bounding boxes of the partitioned geometry elements, and $\mathrm{diam}(\cdot)$ and $\mathrm{dist}(\cdot,\cdot)$ denote the Euclidian diameter and the minimum distance between these bounding boxes, respectively. The parameter $\eta$ is a real positive number that depends on the problem and is chosen empirically.

To apply H-matrices to the CBIE formulation, we begin by constructing the hierarchical block cluster-tree through recursive geometry partitioning based on the Euclidean distances between pairs of observation and source (quadrature) points by testing the admissibility criteria (21) at each step. The blocks that do not satisfy (21) are recursively partitioned until the number of basis/testing functions becomes less than the prescribed leaf-size $n_{\mathrm{leaf}}$, where $n_{\mathrm{leaf}}$ is an empirically chosen parameter that controls the tree depth. Blocks that satisfy (21) are labelled as admissible and the rest are labelled as inadmissible blocks. The matrix entries within inadmissible blocks are stored in a full-rank matrix format without any approximation. Admissible blocks are stored in a low-rank approximate form, computed for example using the Adaptive Cross Approximation algorithm [22].

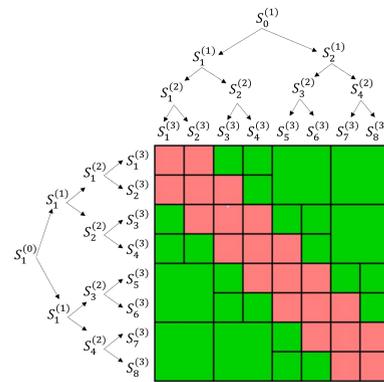

**FIGURE 2.** A 3-level partitioned representation of the H-matrix for the interaction of CBIE source and observation points on a surface.

To demonstrate the hierarchical matrix representation of the NM-CBIE matrix, a block cluster tree representation of the interaction matrix of the observation and source (quadrature) points on the surface is given in Figure 2 for a three-level tree depth. In Fig. 2, green blocks are admissible, while red blocks are inadmissible. After creating the H-matrices, the system can be solved using fast direct H-matrix block H-LU decomposition followed by a block forward and backward substitution applied to the H-matrix structure [23].

In the next two subsections we consider and describe two different approaches for applying the H-matrix framework to the CBIE, each of which is based on a different ordering of the matrix equations shown in Figure 3. In the figure, both a component-wise order (Figure 3a) and a repeated points order (Figure 3b) are shown. Due to its simplicity and improved factorization timings, we use the "entire matrix implementation" approach for all the numerical results in this work.

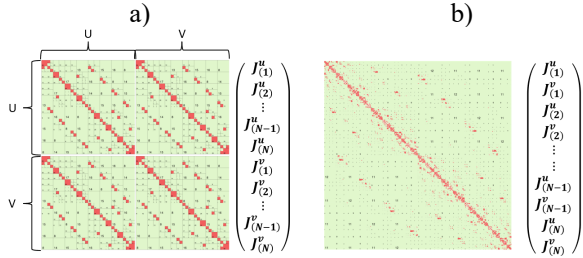

**FIGURE 3.** H-matrix structure obtained for the compression approaches considered herein: a) An individual H-matrix applied to each of the 4 component interaction blocks b) H-matrix applied to the entire set of interactions.

### A. COMPONENT-WISE IMPLEMENTATION

The component-wise approach applies a unique H-matrix representation to each of the component-to-component interaction blocks. As an example, consider the MFIE for a PEC object, where we solve the second row of (10), resulting in evaluating (19) when **M=0**. This leads to a 2×2 block matrix (Fig. 3a)

$$\begin{bmatrix} \mathbf{K}_{uu}^{J,H} & \mathbf{K}_{uv}^{J,H} \\ \mathbf{K}_{vu}^{J,H} & \mathbf{K}_{vv}^{J,H} \end{bmatrix} \begin{bmatrix} \mathbf{J}^u \\ \mathbf{J}^v \end{bmatrix} = \begin{bmatrix} \mathbf{b}_1 \\ \mathbf{b}_2 \end{bmatrix}, \quad (22)$$

where $\mathbf{b}_1 = [-\mu_1 \mathbf{a}_v^p \cdot \mathbf{H}^{p,\text{inc}}]^T$ and $\mathbf{b}_2 = [\mu_1 \mathbf{a}_u^p \cdot \mathbf{H}^{p,\text{inc}}]^T$. With the component-wise approach, we associate an H-matrix with each of the four matrix blocks. As a whole, the overall LU decomposition of (22) can be represented as a block decomposition,

$$\mathbf{K}^{\mathrm{H}} = \begin{bmatrix} \mathbf{K}_{11}^{\mathrm{H}} & \mathbf{K}_{12}^{\mathrm{H}} \\ \mathbf{K}_{21}^{\mathrm{H}} & \mathbf{K}_{22}^{\mathrm{H}} \end{bmatrix} = \begin{bmatrix} L_{11} & 0 \\ L_{21} & L_{22} \end{bmatrix} \cdot \begin{bmatrix} U_{11} & U_{21} \\ 0 & U_{22} \end{bmatrix} = L \cdot U, \quad (23)$$

where $\mathbf{K}_{ij}^{\mathbf{H}}$ are the H-matrix representations of its subblocks. In this approach the H-matrix operations (LU decomposition and back substitution) are applied to each subblock separately and then assembled into the $\mathbf{K}^{\mathbf{H}}$. For example, $\mathbf{K}_{11}^{\mathbf{H}} = L_{11} U_{11}$ is the H-LU decomposition of $\mathbf{K}_{11}^{\mathbf{H}}$. To solve for the surface current densities of (22) using the constructed H-matrix system, the solution of the lower triangle system is first obtained by applying H-matrix forward substitution to the following system

$$L \begin{bmatrix} X_1 \\ X_2 \end{bmatrix} = \begin{bmatrix} b_1 \\ b_2 \end{bmatrix} \quad (24)$$

and then the solution of (24) is used to obtain the unknown densities by applying H-matrix back-substitution as

$$U \begin{bmatrix} \mathbf{J}^u \\ \mathbf{J}^v \end{bmatrix} = \begin{bmatrix} X_1 \\ X_2 \end{bmatrix}. \quad (25)$$

The subblocks $L_{ij}$ and $U_{ij}$ in (23) are computed through recursive calls of block H-LU [22]. The first unknown $X_1$ in (24) can then be computed from

$$L_{11} X_1 = b_1, \quad (26)$$

and the second unknown $X_2$ is obtained by solving

$$L_{22} X_2 = b_2 - L_{21} X_1, \quad (27)$$

in both cases using H-matrix based forward substitution. Finally, the unknown surface current densities in (25) are obtained by solving

$$\begin{bmatrix} U_{11} & U_{12} \\ & U_{22} \end{bmatrix} \begin{bmatrix} \mathbf{J}^u \\ \mathbf{J}^v \end{bmatrix} = \begin{bmatrix} X_1 \\ X_2 \end{bmatrix} \quad (28)$$

in a similar manner using H-matrix based back substitution.

### B. ENTIRE MATRIX IMPLEMENTATION

The CBIE observation and source points result in four scalar interaction matrices ($uu, uv, vu, vv$) for the MFIE or sixteen scalar interactions as shown in equation (19) for the N-Müller formulation. Rather than represent each of these interaction sets with their own H-matrix, the clustering can be applied once to the entire system matrix by forming a single cluster-tree that treats each vector component as an individual point, where points with the same physical location are treated as repeated points (two for MFIE and four for N-Müller). H-matrix operations are thus applied to the entire matrix treating it as a scalar problem without considering the underlying block-wise structure. In the case of the MFIE, the points are repeated twice to account for $\mathbf{J}^u$ and $\mathbf{J}^v$ and since the unknowns are ordered with the u and v components interleaved, the resulting matrix has a diagonally dominant structure (Fig. 3b).

## IV. NUMERICAL RESULTS

To validate the proposed H-matrix accelerated CBIE formulation we present results for both PEC scattering based on the MFIE and dielectric scattering using the N-Müller formulation. All simulations were performed on a server with dual AMD 7763 CPUs (128 cores) and 1TB of RAM.

### A. MFIE

As a first example, we present scattering from a PEC sphere centred at the origin with radius $3\lambda$ illuminated by a plane wave propagating in the $+\hat{z}$ direction and polarized in the $+\hat{x}$ direction with wavelength $\lambda = 1$m and unit amplitude. The number of patches is varied from 24 to 4056 with 10x10 discretization points per patch; the number of unknowns ranges from 4800 to 811,200. A maximum H-matrix Compression Rate (CR) of 98.8% (measured as one minus the ratio of the required H-matrix memory to memory that would be required to store an uncompressed version of the system) was achieved for the case with 811,200 unknowns. The performance of our proposed accelerated solver in terms of time and memory consumption compared against the unaccelerated dense matrix case is shown in Fig. 4. The CPU time required to fill the H-Matrix and the corresponding full dense matrix scales as $O(N\log N)$ and $O(N^2)$ respectively. The time to factor the H-Matrix and the corresponding dense matrix approximately scales as $O(N\log^2 N)$ and $O(N^3)$ respectively. Although the factorization time for the H-Matrix deviates from $O(N\log^2 N)$ scaling (likely due to the ranks of the off-diagonal blocks growing faster for higher frequency problems), significant reduction of CPU time can still be observed in comparison with the factorization time of the $O(N^3)$ dense matrix scaling. The time taken for backsubstitution after the LU factorization has completed scales as $O(N\log N)$. Finally, Fig. 4(b) shows that the memory consumption for the H-Matrix case scales approximately as $O(N\log N)$, whereas that of the corresponding dense matrix case scales as $O(N^2)$.

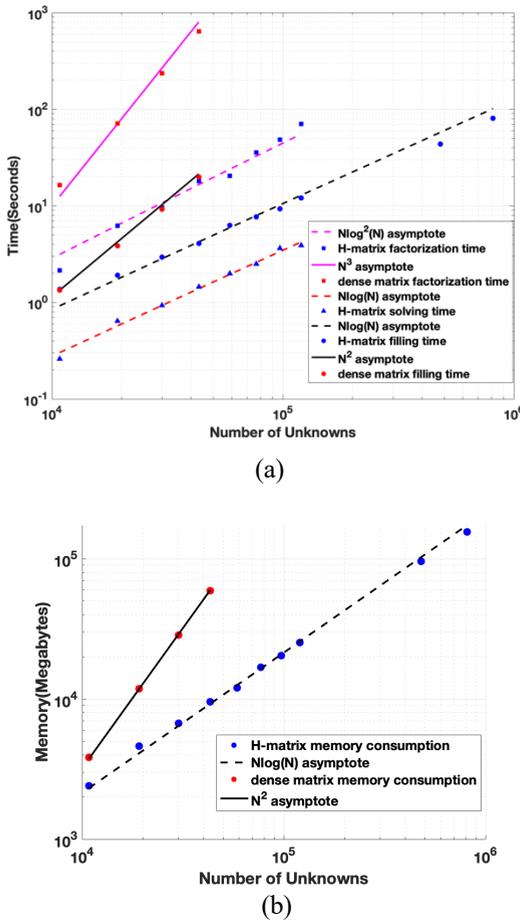

(a)

(b)

**FIGURE 4.** Performance comparison of H-matrix direct solver vs. dense matrix direct solver. (a) CPU time complexity. (b) Memory consumption.

**TABLE 1.** Error of the density solved by iterative solver and proposed direct solver for different number of unknowns.

| # of unknowns | 43200 | 76800 | 120000 | 172800 |
|---|---|---|---|---|
| Iterative | 1.91e-01 | 1.47e-02 | 7.39e-03 | 4.00e-03 |
| Proposed | 1.91e-01 | 1.41e-02 | 6.69e-03 | 2.13e-03 |

For the next example, we compared the accuracy of the solved densities from scattering by $20\lambda$ diameter PEC sphere between the uncompressed GMRES iterative solver and the proposed H-matrix accelerated direct solver. The tolerance for the ACA compression algorithm to build the H-matrix and for the GMRES convergence was set to 1e-4 in both cases. Table 1 shows that the accuracy compared to the analytical Mie series solution is very similar in both cases with the total number of unknowns ranging from 43,200 up to 172,800, indicating that the H-matrix direct solver remains accurate even for electrically large geometries.

To demonstrate the performance advantages of the proposed CBIE method over the conventional LCN method, we solved a plane wave scattering problem with a $1\lambda$-diameter PEC sphere using both methods with the same exact mesh, polynomial order, and number of unknowns. The sphere was partitioned into 384 patches with 4x4 discretization points per patch. It took 362 seconds to build the H-matrix with the conventional LCN method and only 4.6 seconds to build it using the CBIE method for this example, corresponding to an 83x speedup. As with all the other examples, both solvers were evaluated on the same server with dual AMD 7763 CPUs (128 cores). This significant reduction in time consumption clearly demonstrates the benefits of using the CBIE method to build the H-matrix system.

In order to show the capability of our proposed solver for handling objects with more complex geometrical features, we consider scattering of a plane wave from a B2 aircraft model with $16\lambda$ wingspan. Fig. 5(a) shows the propagation direction $\hat{k}$ of the incident field and the real part of the x-component of the induced current density **J** on the surface. Fig. 5(b) plots the radar cross-section (RCS) versus $\theta$ for the $\varphi = 0°$ plane, computed both by using the proposed solver as well as an unaccelerated iterative solver for comparison using 51,600 unknowns. The H-matrix achieved a 88.2% CR and the two solutions match closely. Finally, we consider excitation of a $2\lambda$ diameter metallic spiral antenna by a unit amplitude x-polarized dipole source located at the origin. The spiral was discretized with 11,520 unknowns and achieved a 72% CR. Fig. 6 depicts the magnitude of the solved current density **J** on the antenna's surface.

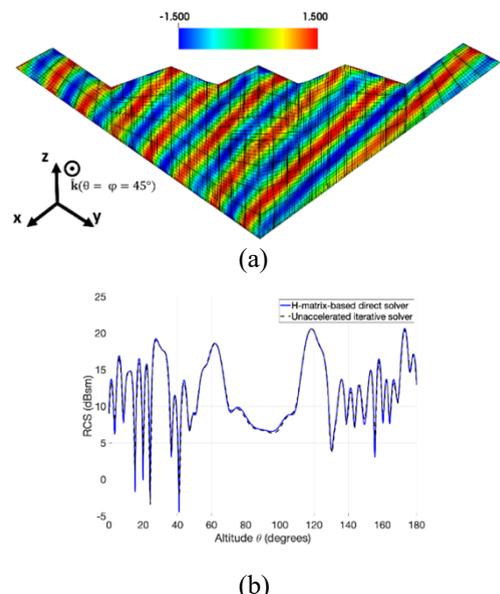

(a)

(b)

**FIGURE 5.** (a) Re($J_x$) on the surface of the B2 aircraft. (b) RCS for $\varphi=0°$ comparing H-Matrix direct solver vs. unaccelerted iterative solver.

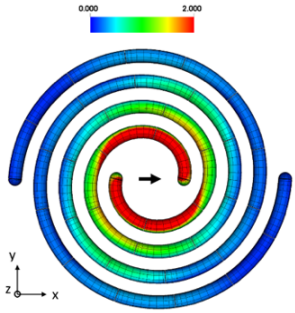

**FIGURE 6.** The magnitude of the current density |J| on the surface of the antenna.

### B. N-Müller Formulation

In this section, scattering from two different dielectric objects is computed to illustrate the effectiveness of our proposed solver when applied to the Müller formulation. The first example we consider is a dielectric sphere in free-space centred at the origin with diameter $4\lambda_e$ illuminated by a plane wave propagating in $+\hat{z}$ direction and polarized in $+\hat{x}$ direction with wavelength of the exterior free-space region $\lambda_e = 1$ m and unit amplitude. The surface of the sphere is partitioned into 96 patches with 10x10 discretization points per patch, resulting in 38,400 unknowns. It is shown in [24] that sharp resonances occur for a dielectric sphere as its size or refractive index is varied, which, in our example, is realized by sweeping the permittivity of the sphere $\varepsilon_d$ around the value corresponding to the resonance. The values of $\varepsilon_d$ corresponding to physical resonances can be found according to the methodology in [24]. As shown in Fig. 7(a), the scattering is computed from $\varepsilon_d = 10.3$ to $\varepsilon_d = 10.7$ solved both iteratively using GMRES and directly using our proposed H-matrix accelerated solver, and the observed resonance matches closely with the predicted $\varepsilon_d = 10.5082$. The time required for a direct solution stays stable (270 seconds) even if the permittivity is extremely close to the resonance ($\varepsilon_d = 10.5074$), while the iterative solver struggles to converge, demonstrating performance that is highly dependent on the resonance. Fig. 7(b) shows the current density distribution corresponding to $\varepsilon_d = 10.5$, and Table 2 compares the errors of the iterative and direct solvers against the analytical Mie series solution for four labeled points in Fig. 7(a), highlighting that the accuracy remains consistent even near the resonance.

To demonstrate the effectiveness of the approach for problems with high permittivity contrast, we consider a second example of a dielectric sphere in free-space with $1\lambda_e$ diameter and permittivities ranging from 50 to 300. The same excitation is used as the first example. Fig. 8 shows the total time required to solve the problem using both the iterative GMRES solver on the uncompressed system and our proposed H-matrix accelerated direct solver. The accuracy of the solved density was controlled to be better than 1e-2 in all cases by setting the ACA tolerance to 1e-8 for the H-matrix solver and the GMRES tolerance to 1e-4 for the iterative solver. It can be seen that the number of iterations and the solution time grows rapidly for the iterative solver as the permittivity increases, while the time needed by the direct solver remains nearly constant even for permittivity values as high as 300. This highlights the ability of the direct solver to handle challenging problems with extremely high permittivities efficiently.

The last example is a dielectric humanoid bunny CAD model generated by a commercial software package [25]. Due to the complexity of this geometry, it can be challenging to iteratively solve when the permittivity is high. Fig. 9(a) compares the time required by the H-matrix based direct solver and the iterative solver. As the permittivity increases, the time cost for the iterative solver increases rapidly due to requiring larger numbers of iterations to reach convergence, while the time to solution of the direct solver does not grow appreciably even under the condition of maintaining the same accuracy in both solutions. Fig. 9(b) shows the current density distribution corresponding to $\varepsilon_d = 6$. The CPU time and memory complexity are also investigated for this example. The number of discretization points per patch, corresponding to the basis order used, is varied from 11x11 to 25x25 with 402 patches, which results in the number of unknowns ranging from 194,568 to 1,005,000. A maximum H-matrix Compression Rate (CR) of 98.5% was achieved for the case with 1,005,000 unknowns. Fig. 10(a) shows that the H-matrix building and factorization time scales as $O(N\log N)$ and $O(N\log^2 N)$ approximately. Fig. 10(b) shows that both the H-matrix building memory and total memory (building and factorization) scales as $O(N \log N)$.

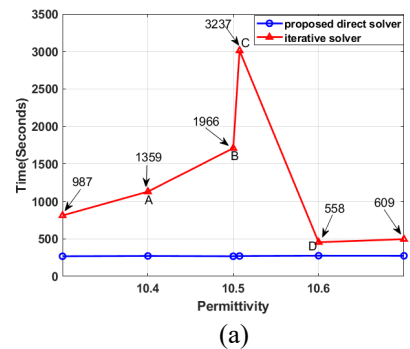

(a)

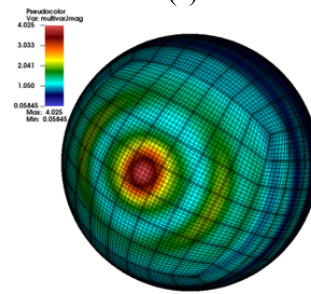

(b)

**FIGURE 7.** a) Time consumption of the iterative solver versus the proposed H-matrix accelerated direct solver for various permittivity values close to resonance. The number of iterations is also shown for each data point b) Surface J distribution on a dielectric sphere with $4\lambda_e$ diameter and permittivity $\varepsilon = 10.5$.

**TABLE 2.** Error of the density solved by iterative solver and proposed direct solver for data points A, B, C and D shown in Figure 7(a).

|  | A | B | C | D |
|---|---|---|---|---|
| iterative | 1.26e-02 | 3.39e-02 | 8.70e-02 | 3.64e-03 |
| proposed | 1.26e-02 | 3.39e-02 | 8.69e-02 | 3.64e-03 |

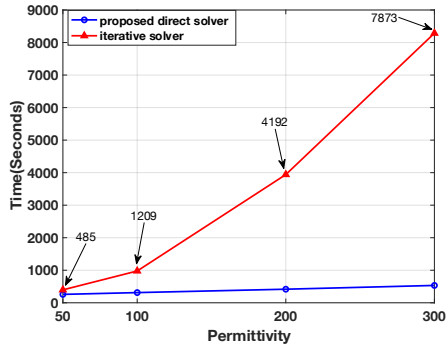

**FIGURE 8.** Time consumption of the iterative solver versus the proposed H-matrix accelerated direct solver for high permittivity values ranging from 50 to 300. The number of iterations is also shown for each data point.

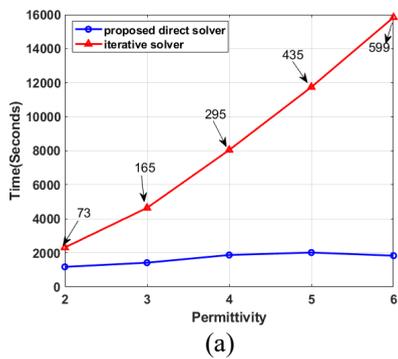

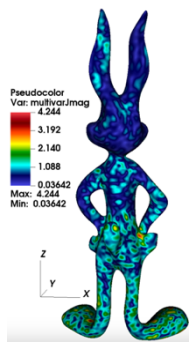

**FIGURE 9.** a) Surface J distribution on a 16λ CAD humanoid bunny model with permittivity $\varepsilon = 6$. b) Time consumption of the iterative solver versus the proposed H-matrix accelerated direct solver for various permittivity values. The number of iterations is also shown for each data point.

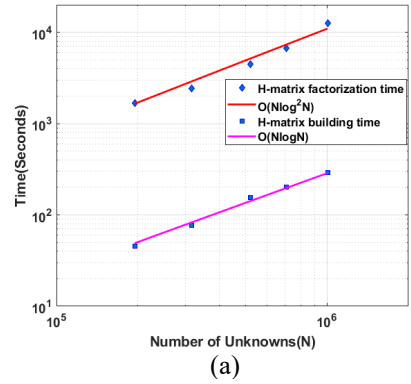

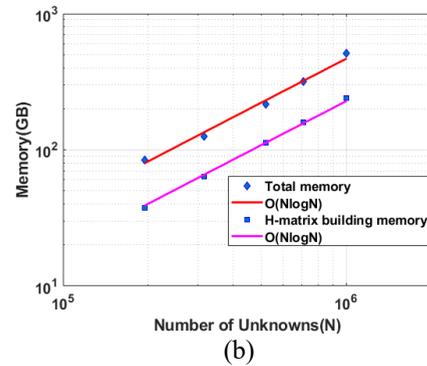

**FIGURE 10.** a) CPU time b) memory complexity of the H-matrix accelerated direct solver for solving the dielectric bugs bunny model with various number of unknowns.

## V. CONCLUSION

An accelerated direct solver for the recently introduced CBIE method was presented using H-matrices for both MFIE and N-Müller formulations. Results were reported on the performance of the approach in terms of matrix fill, factorization, solve time, and memory requirements, demonstrating a new computational framework, which is both computationally efficient and resilient to problems that are challenging for iterative approaches, such as scenarios with physical resonances or high permittivities. The matrix fill time is fast due to the use of the CBIE method, requiring only seconds to minutes to complete even for problems approaching 1 million unknowns and achieving more than 86x speedup over a conventional LCN approach. Scattering from a NURBS-based airplane CAD model, a metallic spiral antenna, and a dielectric bunny, also demonstrated the formulation's capabilities for handling complex geometries. Although all the simulations done in this paper were parallelized on a single machine with shared memory using OpenMP, it should be straightforward to efficiently parallelize the H-matrix building part of the solver on a multi-node HPC cluster using MPI since there is very little data communication involved at this stage. However, parallelizing the LU-decomposition part across a multi-node cluster is tougher due to the larger amount of communication required between nodes and challenges in efficiently load balancing the worker threads. We expect that task programming libraries such as StarPU can be used to solve some of these issues due to their capability for optimizing data communication and task scheduling among different nodes. Extending this work to scale efficiently on large-scale HPC clusters or alternative hardware architectures such as GPUs is the subject of on-going and future work.

## ACKNOWLEDGMENT

The authors are grateful for support by the Air Force Office of Scientific Research (FA9550-20-1-0087) and by the National Science Foundation (CCF-1849965, CCF-2047433).